\documentclass[12pt]{article}
\usepackage{amsmath,amsthm,amssymb}

\renewcommand{\theequation}{\mbox{\arabic{section}.\arabic{equation}}}

\newtheorem{theorem}{Theorem}[section] 
\newtheorem{lemma}{Lemma}[section]
\newtheorem{cor}{Corollary}[section]

\newtheorem{example}{Example}[section]

\newtheorem{proposition}{Proposition}[section]

\numberwithin{equation}{section}

\newcommand{\D}{{\rm d}}
\newcommand{\dx}{\, \D x}

\newcommand{\ds}{\, \D s}
\newcommand{\dr}{\, \D r}
\newcommand{\dt}{\, \D t}

\newcommand{\dis}{\displaystyle}

\newcommand{\msp}{\;\;}
\newcommand{\fsp}{\quad\;}
\newcommand{\psp}{\,}

\newcommand{\arsp}{\hspace*{\arraycolsep}}

\newcommand{\rz}{\mathbb{R}}
\newcommand{\nz}{\mathbb{N}}

\newcommand{\eps}{\varepsilon}
\newcommand{\iom}{\int_{\Omega}}

\newcommand{\klauf}{\left(\begin{array}}
\newcommand{\klzu}{\end{array}\right)}

\newcommand{\udel}{u_{\delta}}
\newcommand{\gam}[2]{\Gamma_{#1,\delta}^{#2}}
\title{Splitting-type variational problems with linear growth conditions}
\author{Michael Bildhauer \& Martin Fuchs}
\date{}

\newcommand{\reff}[1]{(\ref{#1})}

\usepackage{hyperref}

\hypersetup{
   pdfnewwindow=true,
   colorlinks=false,
   linkcolor=black,
   citecolor=green,
   filecolor=magenta,
   urlcolor=cyan,
}

\begin{document}

\parindent0em
\maketitle

\newcommand{\op}[1]{\operatorname{#1}}
\newcommand{\bv}{\op{BV}}
\newcommand{\mub}{\overline{\mu}}

\newcommand{\hypref}[2]{\hyperref[#2]{#1 \ref*{#2}}}
\newcommand{\hypreff}[1]{\hyperref[#1]{(\ref*{#1})}}

\newcommand{\ob}[1]{^{(#1)}}

\begin{abstract}
Regularity properties of solutions to variational problems are established for a broad class 
of strictly convex splitting-type energy densities of the principal form $f$: $\rz^2 \to \rz$,
\[
f(\xi_1,\xi_2) = f_1\big( \xi_1 \big) + f_2\big( \xi_2 \big) \psp ,
\]
with linear growth. As a main result it is shown that, regardless of a corresponding property of $f_2$,
the assumption ($t\in \rz$)
\[
c_1 (1+|t|)^{-\mu_1} \leq f_1''(t) \leq c_2\psp ,\fsp 1 < \mu_1 < 2\psp ,
\]
is sufficient to obtain higher integrability of $\partial_1 u$ for any finite exponent. 
We also include a series of variants of our main theorem. We finally note that 
similar results in the case $f$: $\rz^n \to \rz$
hold with the obvious changes in notation.\footnote{AMS-Classification: 49J45, 49N60}
\end{abstract}

\parindent0ex
\section{Introduction}\label{intro}

In our paper we discuss variational problems of linear growth with densities
which do not belong to the class of $\mu$-elliptic energies introduced first in
\cite{BF:2003_1}.\\

Guided by linear growth examples of splitting-type, which to our knowledge
are not systematically studied up to now, we are led to quite general hypotheses
which still guarantee some interesting higher regularity properties of generalized solutions.\\

Before going into details, let us fix the framework of our considerations: in what
follows $\Omega$ denotes a bounded Lipschitz domain in $\rz^n$,
$n \geq 2$, and we consider a function $u_0$: $\Omega \to \rz$ such that\footnote{Using a suitable approximation 
(see, e.g., \cite{Bi:2003_2} for more details), it is possible to suppose $u_0\in W^{1,1}(\Omega)\cap L^\infty(\Omega)$.}
\begin{equation}\label{intro 1}
u_0 \in W^{1,2}(\Omega) \cap L^\infty(\Omega) \psp .
\end{equation}

We then are interested in  the variational problem
\begin{equation}\label{intro 2}
J[u] := \iom f(\nabla u) \dx \to \min \fsp\mbox{in}\msp u_0 +W^{1,1}_0(\Omega) 
\end{equation}
for a strictly convex energy density $f$: $\rz^n \to [0,\infty)$ of class $C^2$
satisfying
\begin{equation}\label{intro 3}
a_1 |\xi| - a_2 \leq f(\xi) \leq a_3 |\xi| + a_4 \psp ,\fsp \xi\in \rz^n \psp ,
\end{equation}
with suitable constants $a_1$, $a_3 > 0$, $a_2$, $a_4 \geq 0$.\\ 

Condition \reff{intro 3} causes the well-known problems concerning the existence and the regularity of
solutions to \reff{intro 2}, which means that \reff{intro 2} has to be replaced by a
relaxed variant. For the general framework of this approach we refer, e.g., to the monographs
\cite{Gi:1984_1}, \cite{GMS:1998_1}, \cite{GMS:1998_2}, \cite{AFP:2000_1} and
\cite{Bi:1818}, where the reader will find a lot of further references as well as a definition
of the underlying spaces such as $L^p(\Omega)$, $W^{1,p}(\Omega)$, $\bv(\Omega)$ and their
local variants.\\

Quoting \cite{AFP:2000_1}, Theorem 5.47, the natural extension of \reff{intro 2} reads as
\begin{eqnarray}\label{intro 4}
K[w]&:=& \iom f(\nabla^a w) \dx + \iom f_\infty\Bigg(\frac{\nabla^s w}{|\nabla^s w|}\Bigg) \D |\nabla^s w|
\nonumber\\
&&+ \int_{\partial\Omega} f_\infty \big((u_0-w)\mathcal{N}\big) \D \mathcal{H}^{n-1}
\to \min \fsp\mbox{in}\msp \bv(\Omega) \psp .
\end{eqnarray}
Here $\nabla w = \nabla^a w \mathcal{L}^n + \nabla^s w$ is the Lebesgue decomposition
of the vector measure $\nabla w$ with respect to the $n$-dimensional Lebesgue measure
$\mathcal{L}^n$, $f_\infty$ is the recession fuinction of $f$, i.e.
\[
f_\infty(\xi) := \lim_{t\to \infty} \frac{1}{t}f(t\xi)\psp ,\fsp \xi \in \rz^n\psp ,
\]
$\mathcal{H}^{n-1}$ is Hausdorff's measure of dimension $n-1$ and $\mathcal{N}$
denotes the outward unit normal to $\partial \Omega$.\\

We summarize some important results concerning the relations between problems \reff{intro 2}
and \reff{intro 4} in the following proposition (compare, e.g., 
the pioneering work \cite{GMS:1979_1} and \cite{GMS:1979_2} in the minimal surface case,
and, e.g., the papers \cite{Se:1985_1}, \cite{Se:1996_1}, where the mechanical point of view
is discussed by introducing the stress tensor as the unique solution of the dual problem).

\begin{proposition}\label{intro prop 1}
Let \reff{intro 1} and \reff{intro 3} hold.
\begin{enumerate}
\item Problem \reff{intro 4} admits at least one solution $u \in \bv(\Omega)$.

\item It holds
\[
\inf_{u_0 + W^{1,1}_0(\Omega)} J = \inf_{\bv(\Omega)} K \psp .
\]

\item We have the following characterization:
\[
\begin{array}{l}
\mbox{$u \in \bv(\Omega)$ is $K$-minimizing}\arsp  \Leftrightarrow\\[1ex]
u\in \mathcal{M} :=
\Big\{v \in L^1(\Omega):\psp \mbox{$v$ is a $L^1(\Omega)$-cluster point
of some}\\[1ex] 
\qquad\qquad\qquad\mbox{$J$-minimizing sequence from $u_0 + W^{1,1}_0(\Omega)$}\Big\}\psp .
\end{array}
\]
\end{enumerate}

\end{proposition}

Since Proposition \ref{intro prop 1} in particular guarantees the existence of generalized
solutions to problem \reff{intro 2}, i.e.~of functions $u\in \bv(\Omega)$ solving \reff{intro 4},
one may ask for their regularity properties.\\

Here a variety of results is available concerning densities $f$ of linear growth such that
we have in addition
\begin{equation}\label{intro 5}
c_1 \big(1+|\xi|\big)^{-\mu} |\eta|^2 \leq D^2f(\xi) (\eta,\eta) \leq c_2 \big(1+|\xi|\big)^{-1}
|\eta|^2
\end{equation}
with exponent $\mu >1$ and for constants $c_1$, $c_2 > 0$.\\

Condition \reff{intro 5} is known as the $\mu$-ellipticity property of $f$.
Assuming at least $\mu \leq 3$, the reader will find regularity results e.g.~in the paper
\cite{BF:2003_1}, in \cite{Bi:2002_1} the case of bounded solutions is covered. We also mention the work of Marcellini and Papi \cite{MP:2006_1}
and the paper \cite{BS:2013_1} together with the references quoted therein.\\

Condition \reff{intro 5} is mainly motivated by 
\begin{itemize}
\item energy densities $f$ of minimal surface type, i.e.
\begin{equation}\label{intro 6}
f(\xi) := \big( 1+|\xi|^k\big)^{\frac{1}{k}} \psp ,\fsp k > 1\psp ,
\end{equation}
\item or by densities of the form $f(\xi) = \Phi_\mu\big(|\xi|\big)$, where for $\mu > 1$ we let ($t\geq 0$)
\begin{eqnarray}\label{intro 7}
\Phi_\mu(t) &:=& (\mu-1) \int_0^t\int_0^s (1+r)^{-\mu}\dr\ds\psp , \nonumber\\[1ex]
& =& \left\{\begin{array}{rcl}
\dis t - \frac{1}{2-\mu} (1 + t)^{2- \mu} +  \frac{1}{2- \mu}&\mbox{if}& \mu\not=2\psp ,\\[3ex]
\dis t - \ln(1+t)&\mbox{if}& \mu =2 \psp .
\end{array}\right.
\end{eqnarray}
\end{itemize}

Note that recent strain-limiting elastic models with linear growth are strongly related to the class given in \reff{intro 6}
(see, for instance, \cite{BBMS:2017_1}, \cite{BMS:2015_1}, \cite{BMRS:2014_1} and \cite{BMRW:2015_1}).\\

Let us have a closer look at the second kind of examples.
We carefully have to distinguish between the functions $f(\xi) = \Phi_\mu\big(|\xi|\big)$ 
defined on $\rz^n$ and the functions $\Phi_\mu(t)$ depending on one variable in the sense that
we have with optimal exponents the inequalities \reff{intro 5} for $f$, whereas for all $t\geq 0$
\begin{equation}\label{intro 8}
c_1 (1+t)^{-\mu} \leq \Phi''_{\mu}(t) \leq c_2 (1+t)^{-\mu} \psp ,\fsp \mu > 1\psp ,
\end{equation}
$c_1$, $c_2 >0$. Note that both the exponent occuring in the upper bound of \reff{intro 5} and the one of \reff{intro 8}
are relevant quantities entering the regularity proofs in an essential way.\\

For instance, in the recent paper \cite{BF:2020_2}, the authors benefit from the radial structure
of a solution which roughly speaking means that the general ellipticity condition \reff{intro 5}
can be replaced by the estimate on the right-hand side of \reff{intro 8}.\\

However, without using the radial structure of a particular solution, it is no longer possible 
to benefit that much from the right-hand side of \reff{intro 8}.\\

This becomes evident with the following auxiliary lemma which, roughly speaking, states that the right-hand side of
\reff{intro 5} with exponent $-1$ gives the best possible estimate. 
Without reducing the problem by, e.g., symmetry properties of the solution, 
the balancing condition \reff{intro 8} may not serve as an additional tool for proving the regularity of solutions.

\begin{lemma}\label{intro lem 1}
Let $n \geq 2$ and consider a density $f$ of class $C^2$ such that \reff{intro 3} and
\begin{equation}
c_1\big(1+|\xi|\big)^{-\mu} |\eta|^2 \leq D^2f(\xi)(\eta,\eta) \leq c_2
\big(1+|\xi|\big)^{-\kappa} |\eta|^2\label{intro 9}
\end{equation}
hold for all $\xi$, $\eta \in \rz^n$ with constants $c$, $\overline{c} >0$,
and with exponents $\mu >1$, $\kappa \leq \mu$.
Then we have
\[
\kappa \leq 1\psp .
\]
\end{lemma}
The proof of Lemma \ref{intro lem 1} is postponed to the appendix.\\

Once that $\kappa = 1$ is seen to be the best possible choice in \reff{intro 9}, the question
arises, whether this yields a sufficiently broad class of examples. However, this is not the case if we 
like to include some kind of splitting structure in our considerations:

\begin{example}\label{intro ex 1}
For the sake of simplicity let $n=2$ and consider the energy density of splitting-type
\begin{equation}\label{intro 10}
f(\xi_1,\xi_2) = \Phi_{\mu_1}\big(|\xi_1|\big) + \Phi_{\mu_2} \big(|\xi_2|\big) \psp ,
\fsp \mu_1 ,\psp \mu_2 > 1\psp .
\end{equation}
Then  we merely have ($\mu := \max \{\mu_1,\mu_2\}$, $c_1$, $c_2>0$)
\begin{equation}\label{intro 11}
c_1 \big( 1+|\xi|\big)^{-\mu} |\eta|^2 \leq D^2 f(\xi)(\eta,\eta) \leq
c_2 |\eta|^2\psp ,\fsp \xi, \psp\eta \in \rz^2\psp ,
\end{equation}
and the estimate on the r.h.s.~can not be improved
as we recognize in the case $|\xi_1|\to \infty$ together with $\xi_2=const$, $\eta_1 = 0$.\\
\end{example}

We emphasize that for this kind of splitting examples the balancing condition \reff{intro 8}
can not be exploited for improving the condition \reff{intro 11}, hence instead of $\Phi_\mu$
we may as well consider functions $\psi_\mu$ of the type (adapted to \reff{intro 11})
 
\begin{equation}\label{intro 12}
c_1 (1+|t|)^{-\mu} \leq \psi''_{\mu}(t) \leq c_2  \psp ,\fsp \mu > 1\psp , \fsp t\in \rz\psp ,
\end{equation}
with constants $c_1$, $c_2>0$ still having linear growth.
We like to mention that we also do not rely on a $\Delta_2$-condition similar to, e.g., 
(1.8) of \cite{BBM:2018_1}.\\

One may ask whether this generalization is a kind of artefact without providing new relevant examples.
The examples sketched in the appendix illustrate that this is not the case.\\

With \reff{intro 11} and \reff{intro 12} we are lead to our main theorem on the regularity of
solutions.

\begin{theorem}\label{intro theo 1}
Assume \reff{intro 1}, let $n=2$ and consider a density $f$ of the form 
\begin{equation}\label{intro 13}
f(\xi_1,\xi_2) = f_{1}\big( \xi_1 \big) + f_{2} \big( \xi_2 \big) 
\end{equation}

with strictly convex functions $f_1$, $f_2$ satisfying \reff{intro 3}.
Moreover, suppose that there exist exponents $\mu_i >  1$ such that for $i=1$, $2$
\begin{equation}\label{intro 14}
c_i \big(1+|t|\big)^{-\mu_i} \leq f_i''(t) \leq \bar{c}_i\psp , \fsp t \in \rz \psp , 
\end{equation}
holds with constants $c_i$, $\bar{c}_i > 0$.\\

Let $\mu_1 < 2$. Then there exists a generalized minimizer $u \in \mathcal{M}$ such that
\[
\partial_1 u \in L^\chi_{\op{loc}}(\Omega)\fsp\mbox{for any finite $\chi$}\psp .
\]

With obvious changes in notation, similar results hold in the case $n\geq 3$, e.g.~for the density
\[
f(\xi_1,\dots ,\xi_n) = \sum_{i=1}^n f_i(\xi_i)\psp .
\]
\end{theorem}

On one hand, Theorem \ref{intro theo 1} states that higher integrabilty w.r.t.~a particular direction in the splitting case holds 
provided that the corresponding part of the energy satisfies a sufficient ellipticity condition. No further restriction w.r.t.~the
second direction is imposed.
If we suppose in addition $\mu_2 < 2$, then we expect regular solutions being unique up to a constant. This statement is formulated in
Corollary \ref{intro cor 1} below.\\

On the other hand, condition \reff{intro 11} (with $\mu := \max\{\mu_1,\mu_2\}$ in the splitting case) is very much in the spirit 
of the ellipticity condition 
\begin{equation}\label{intro 15}
c_1 \big( 1+|\xi|\big)^{p-2} |\eta|^2 \leq D^2 f(\xi)(\eta,\eta) \leq
c_2 \big( 1+|\xi|\big)^{q-2} |\eta|^2\psp ,\fsp \xi, \psp\eta \in \rz^2 \psp ,
\end{equation}
for variational problems with anisotropic superlinear growth, where we have
the formal correspondence $q=2$ and $\mu =2-p$.\\

Motivated by the famous counterexample of Giaquinta \cite{Gi:1987_1} there are a lot of contributions to the
regularity theory of solutions which, due to the counterexample, have to impose a suitable relation between the exponents $p$
and $q$. Let us just mention the classical paper \cite{Ma:1989_1}, the reference \cite{ELM:1999_1} on higher integrability
or the recent paper \cite{BM:2019_1}.\\

Note that a series of papers is devoted to the splitting case, for instance 
\cite{Gi:1987_1}, \cite{FS:1993_1}, \cite{Br:2010_1}, \cite{BFZ:2007_1}.\\

In the case of bounded solutions, one suitable relation between $p$ and $q$ implying the regularity of solutions reads as
\[
q < p+ 2
\]
which exactly corresponds to $q=2$, i.e.~$\kappa =0$ and $\mu = 2- p < 2$ in \reff{intro 16} below.\\

In this spirit we have Corollary \ref{intro cor 2} presenting a uniqueness and regularity result for generalized solutions $u \in \mathcal{M}$
in the case that the density $f$ is of linear growth, but not necessarily with splitting structure. The appropriate version of \reff{intro 15}
is stated in \reff{intro 16} and the reader should note that this ellipticity condition does not require the boundedness of
$D^2f(\xi)$. We also like to remark that -- due to the missing splitting structure -- Corollary \ref{intro cor 2} is not a formal 
consequence of Theorem \ref{intro theo 1} or Corollary \ref{intro cor 1} but can immediately deduced by going through the arguments 
presented there.\\

To continue with the main line of splitting-type variational integrals we take Corollary \ref{intro cor 2} as a motivation to
eliminate the boundedness of $f_1''$ in \reff{intro 14} still supposing an analogue to \reff{intro 17} for $\mu_1$. This is done in Corollary \ref{intro cor 3}.\\

To include a broader class of energy densities $f_2$ in our considerations, we have to argue with
negative exponents in the Caccioppoli-type inequality. The price we have to pay is an integrability result up a finite exponent $\chi > 2$
presented in Theorem \ref{intro theo 2}. \\

We note that the above reasoning will be carried over to the consideration of
mixed linear-superlinear problems in the forthcoming paper \cite{BF:2020_4}.\\

\begin{cor}\label{intro cor 1} 
Suppose that the assumptions of Theorem \ref{intro theo 1} hold together with
\[
\mu := \max\big\{\mu_1,\mu_2\big\} < 2 \psp .
\]
Then the relaxed problem \reff{intro 4} admits a solution $u\in C^{1,\alpha}(\Omega)$,
$0 < \alpha < 1$. Moreover, this solution is unique up to additive constants.
\end{cor}

\begin{cor}\label{intro cor 2}
Suppose that we have \reff{intro 1}. Let $f$ satisfy \reff{intro 3} together with
\begin{equation}\label{intro 16}
c_1 \big( 1+|\xi|\big)^{-\mu} |\eta|^2 \leq D^2 f(\xi)(\eta,\eta) \leq
c_2 \big( 1+|\xi|\big)^{-\kappa} |\eta|^2\psp ,\fsp \xi, \psp\eta \in \rz^n \psp ,
\end{equation}
for some constants $c_1$, $c_2 >0$. Suppose further that we have in addition
\begin{equation}\label{intro 17}
\mu < 2 + \kappa
\end{equation}
with exponents $\kappa$, $\mu$ such that
\[
\mu >1 \fsp\mbox{and}\fsp -1 < \kappa \leq 1 \psp 
\]
Then problem \reff{intro 4} admits a solution $u \in \mathcal{M}$ being of class $C^{1,\alpha}(\Omega)$ for any
$0 < \alpha < 1$. Moreover, this solution is unique up to additive constants , which means $v=u+c$ for any
$v\in \mathcal{M}$ with a suitable constant $c\in \rz$.\\
\end{cor}

\begin{cor}\label{intro cor 3}
Theorem \ref{intro theo 1} remains valid if \reff{intro 14} is replaced by the weaker conditions
\begin{eqnarray}\label{vari 1}
c_1 \big(1+|t|\big)^{-\mu_1} \leq &f_1''(t)& \leq \bar{c}_1 \big(1+|t|)^{\varkappa} \psp ,\\[2ex]
c_2 \big(1+|t|\big)^{-\mu_2} \leq &f_2''(t)& \leq \bar{c}_2   \psp , \fsp t \in \rz \psp ,\label{vari 2}
\end{eqnarray}
with constants $c_i$, $\bar{c}_i > 0$, $i=1$, $2$, and with exponents
\[ 
 1 < \mu_1\psp , \fsp 0 \leq \varkappa <  2 - \mu_1 \psp .
\]
\end{cor}
 
 \begin{theorem}\label{intro theo 2}
Suppose that we have the assumptions of Theorem \ref{intro theo 1}, where now \reff{intro 14} is replaced by 
 \begin{eqnarray}\label{vari 3}
c_1 \big(1+|t|\big)^{-\mu_1} \leq &f_1''(t)& \leq \bar{c}_1 \psp , \fsp t \in \rz \psp ,\\[2ex]
c_2 \big(1+|t|\big)^{-\mu_2} \leq &f_2''(t)& \leq \bar{c}_2 \big(1+|t|\big)^\gamma \psp , \fsp t \in \rz \psp .\label{vari 4}
\end{eqnarray}
with constants $c_i$, $\bar{c}_i > 0$, $i=1$, $2$. Moreover we assume that
\begin{equation}\label{vari 5}
0 \leq \gamma < \frac{2-\mu_1}{1+(2-\mu_1)} \psp .
\end{equation}
Then, if $\mu_1 < 2$, there exists a generalized minimizer $u\in \mathcal{M}$ 
\[
\partial_1 u \in L^\chi_{\op{loc}} (\Omega)\fsp\mbox{for some}\msp \chi >2 \psp .
\]
 \end{theorem}

\section{Proof of Theorem \ref{intro theo 1}}\label{proof 1}

We fix some $0 < \delta < 1$ and let
\begin{eqnarray*}
f_{i,\delta}(t) &:= & \frac{\delta}{2} t^2 + f_i(t) \psp , \fsp t\in \rz\psp ,\fsp i=1,\psp 2\psp ,\\[2ex]
f_\delta(\xi) &=& f_{1,\delta}(\xi_1) + f_{2,\delta}(\xi_2) \psp ,\\
&=& \frac{\delta}{2} |\xi|^2 + f(\xi)\psp , \fsp \xi =(\xi_1,\xi_2) \in \rz^2 \psp .
\end{eqnarray*}

We then consider the regularized minimization problem
\renewcommand{\theequation}{\mbox{1.2$_\delta$}}
\begin{equation}\label{delta}
J_\delta[w] := \iom f_\delta(\nabla w) \dx \to \min\fsp\mbox{in}\msp u_0 + W^{1,2}_{0}(\Omega)
\end{equation}\addtocounter{equation}{-1}\renewcommand{\theequation}{\mbox{\arabic{section}.\arabic{equation}}}
with $u_\delta$ denoting the unique solution of \reff{delta}. Following standard arguments
(see \cite{Bi:1818} and a series of well known references quoted therein) one immediately
obtains
\[
\udel \in W^{2,2}_{\op{loc}}(\Omega) \cap C^1(\Omega)
\]
and passing to a subsequence, if necessary, we obtain in the limit $\delta \to 0$
\[
\delta \iom |\nabla \udel|^2 \dx \to 0\psp ,\fsp
\udel \to: u \fsp\mbox{in}\msp L^1(\Omega)\fsp\mbox{with some}\msp u \in \mathcal{M}\psp .
\]
Note using \reff{intro 1} together with the maximum-principle implies
\[
\sup_{\delta} \|\udel\|_{L^{\infty}(\Omega)} < \infty \psp .
\]

We let
\[
\gam{i}{} : = 1+ |\partial_i\udel|^2\psp , \fsp i=1,\psp 2\psp ,  
\]
differentiate the Euler equation
\[
0 = \iom Df_\delta (\nabla \udel) \cdot \nabla \varphi \dx\psp , \fsp \varphi\in C^\infty_0(\Omega)\psp ,
\]
in the sense that we insert $\varphi = \partial_1 \psi$ as test function and obtain for all $\psi \in C^\infty_0(\Omega)$
\begin{equation}\label{proof 1 1}
0 = \iom D^2 f_\delta(\nabla \udel) \big(\nabla \partial_1 \udel, \nabla \psi\big) \dx \psp .
\end{equation}

The first main step in the proof of Theorem \ref{intro theo 1} is to show the following Caccioppoli-type inequality.

\begin{proposition}\label{prop cacc}
Fix $l\in \nz$ and suppose that $\eta \in C^\infty_0(\Omega)$, $0 \leq \eta \leq 1$. Then, given the assumptions of Theorem \ref{intro theo 1},
the inequality
\begin{eqnarray}\label{proof 1 2}
\lefteqn{\iom D^2 f_\delta(\nabla \udel)\big(\nabla \partial_1 \udel, \nabla \partial_1 \udel\big) \eta^{2l} 
\gam{1}{\alpha} \dx}\nonumber \\ 
&& \leq c \iom D^2f_\delta(\nabla \udel) (\nabla \eta,\nabla \eta)\eta^{2l-2} \gam{1}{\alpha +1}\ \dx
\end{eqnarray}
holds for any $\alpha \geq 0$, which in particular implies
\begin{equation}\label{proof 1 3}
\iom \eta^{2l} \gam{1}{\alpha - \frac{\mu_1}{2}} |\partial_{11}\udel|^2 \dx
\leq c \iom |\nabla \eta|^2 \eta^{2l-2}\gam{1}{\alpha +1} \dx \psp .
\end{equation}
Here and in what follows $c=c(\alpha,l)$ denotes a uniform constant independent of $\delta$.
\end{proposition}

\emph{Proof of Proposition \ref{prop cacc}.} We insert the admissible test function
\[
\psi := \eta^{2l} \partial_1 \udel \gam{1}{\alpha}
\]
in \reff{proof 1 1} and obtain
\begin{eqnarray}\label{proof 1 4}
\lefteqn{\iom D^2f_\delta(\nabla \udel) \big(\nabla \partial_1 \udel , \nabla \partial_1 \udel\big)
\eta^{2l} \gam{1}{\alpha}\dx}\nonumber\\
&=& - \iom D^2f_\delta(\nabla \udel)\big(\nabla \partial_1 \udel, \nabla \gam{1}{\alpha}\big)
\partial_1 \udel \eta^{2l}\dx\nonumber\\
&& - \iom D^2f_\delta(\nabla \udel)\big(\nabla \partial_1 \udel, \nabla (\eta^{2l})\big)
\partial_1 \udel \gam{1}{\alpha} \dx  =: S_1+S_2 \psp .
\end{eqnarray}
For $S_1$ we have
\begin{eqnarray*}
S_1 &=& - \alpha \iom D^2f_\delta (\nabla \udel) \big(\nabla \partial_1 \udel, 
\nabla |\partial_1 \udel|^2\big) \gam{1}{\alpha -1} \partial_1\udel \eta^{2l} \dx\\
&=& - 2 \alpha \iom D^2f_\delta (\nabla \udel) \big(\nabla \partial_1 \udel, \nabla \partial_1 \udel\big)
|\partial_1 \udel|^2\gam{1}{\alpha-1}\eta^{2l}\dx \leq 0 \psp ,
\end{eqnarray*}
whenever $\alpha \geq 0$, hence the left-hand side of \reff{proof 1 4} is bounded by $|S_2|$.\\
   
$S_2$ is handled with the help of the Cauchy-Schwarz inequality: for $0 < \eps$ sufficiently small it holds: 
\begin{eqnarray*}
\lefteqn{\iom D^2f_\delta(\nabla \udel) (\nabla \partial_1 \udel, \nabla \eta)\eta^{2l-1} 
\gam{1}{\alpha} \partial_1\udel \dx}\\
&\leq& \eps \iom D^2f_\delta(\nabla \udel) (\nabla \partial_1 \udel ,\nabla \partial_1 \udel) 
\eta^{2l} \gam{1}{\alpha}\dx\\ 
&&+ c(\eps) \iom D^2f_\delta(\nabla \udel) (\nabla \eta,\nabla \eta) 
\eta^{2l-2}\gam{1}{\alpha} |\partial_1 \udel|^2 \dx
\end{eqnarray*}
and absorbing the first term we have shown \reff{proof 1 2}.\\

We now benefit from the splitting structure expressed in \reff{intro 13}, use \reff{intro 14} and estimate the left-hand side
of \reff{proof 1 2} from below
\begin{eqnarray}\label{proof 1 5}
\lefteqn{\iom \eta^{2l} \gam{1}{\alpha - \frac{\mu_1}{2}} |\partial_{11}\udel|^2 \dx}\nonumber\\
&\leq& c \iom  f_{1}''(\partial_1\udel)|\partial_{11}\udel|^2\eta^{2l} \gam{1}{\alpha} \dx 
\leq c \iom  f_{1,\delta}''(\partial_1\udel)|\partial_{11}\udel|^2\eta^{2l} \gam{1}{\alpha} \dx\nonumber\\
&\leq& c \iom \Bigg[ f_{1,\delta}''(\partial_1\udel)|\partial_{11}\udel|^2 
+ f''_{2,\delta}(\partial_2 \udel) |\partial_{12} \udel|^2\Bigg]\eta^{2l} \gam{1}{\alpha}\dx\nonumber\\
&\leq & c \iom D^2f_\delta(\nabla \partial_1\udel, \nabla \partial_1 \udel) \eta^{2l} \gam{1}{\alpha}\dx \psp .
\end{eqnarray}

For the right-hand side of \reff{proof 1 2} we have with $\delta \leq 1$ and recalling \reff{intro 14}
\begin{eqnarray}\label{proof 1 6}
\lefteqn{\iom D^2f_\delta(\nabla \udel)(\nabla \eta,\nabla \eta)\eta^{2l-2} \gam{1}{\alpha +1} \dx}\nonumber\\
&\leq & c \iom \Bigg[f''_{1,\delta}(\partial_1 \udel) |\partial_1 \eta|^2 
+ f''_{2,\delta}(\partial_2 \udel) |\partial_2 \eta|^2\Bigg]\eta^{2l-2} \gam{1}{\alpha +1} \dx\nonumber\\
&\leq & c \iom |\nabla \eta|^2 (1+\delta) \eta^{2l-2}\gam{1}{\alpha +1}\dx \psp .
\end{eqnarray}

By combining \reff{proof 1 2}, \reff{proof 1 5} and \reff{proof 1 6} we have established the claim \reff{proof 1 3}, hence
Proposition \ref{prop cacc}. \hfill\qed\\

Now we are going to discuss the second main ingredient of the proof of Theorem \ref{intro theo 1}.

\begin{proposition}\label{prop high}
Fix some number $\chi >2$ and let
\[
0 < s := \frac{\chi}{2}-1\psp , \fsp \hat{\eps} := 1- \frac{\mu_1}{2} > 0\psp , \fsp \alpha := s - \frac{\hat{\eps}}{2} \psp . 
\]
Then, given the hypotheses of Theorem \ref{intro theo 1}, for $l$ sufficiently large and a local constant $c(\eta,\chi,l)$ we have
\begin{equation}\label{proof 1 7}
\iom \eta^{2l} \gam{1}{s+1} \dx \leq c \Bigg[1+\iom \eta^{2l} \gam{1}{s+\frac{2+\mu_1}{4}} \dx\Bigg] \psp .
\end{equation}
\end{proposition}
\emph{Proof of Proposition \ref{prop high}.} We recall that 
\[
\|\udel\|_{L^\infty(\Omega)} \leq c
\]
with a constant not depending on $\delta$ and estimate
\begin{eqnarray*}
\iom |\partial_1 \udel|^2 \gam{1}{s}\eta^{2l} \dx &= &
\iom \partial_1 \udel \partial_1 \udel \gam{1}{s} \eta^{2l} \dx
= - \iom \udel \partial_1 \Big[\partial_1 \udel \gam{1}{s}\eta^{2l}\Big]\dx\\
&\leq & c \Bigg[ \iom |\partial_{11} \udel| \gam{1}{s} \eta^{2l} \dx +
\iom |\partial_1 \udel| \eta^{2l-1} |\nabla \eta|\gam{1}{s}\dx\\
&&+  \iom \gam{1}{s-1}|\partial_1 \udel|^2 |\partial_{11} \udel| \eta^{2l} \dx \Bigg]\\
&\leq& c \Bigg[ \iom |\partial_{11}\udel| \gam{1}{s} \eta^{2l} \dx \\
&&+ \eps \iom |\partial_1 \udel|^2 \gam{1}{s} \eta^{2l} \dx 
+ c(\eps) \iom |\nabla \eta|^2 \eta^{2l-2}\gam{1}{s} \dx \Bigg]\psp ,
\end{eqnarray*}
where we may choose $\eps >0$ sufficiently small to absorb the second term on the right-hand side.
This means that we have
\begin{equation}\label{proof 1 8}
\iom |\partial_1 \udel|^2 \gam{1}{s} \eta^{2l} \dx \leq
c \Bigg[\iom |\partial_{11} \udel| \gam{1}{s} \eta^{2l} \dx +
\iom |\nabla \eta|^2 \eta^{2l-2}\gam{1}{s}\dx \Bigg] \psp .
\end{equation}

Recalling $\mu_1 < 2$ and using Young's inequality, we estimate the
first term on the right-hand side of \reff{proof 1 8}:
\begin{eqnarray}\label{proof 1 9}
\lefteqn{\iom |\partial_{11} \udel| \gam{1}{s} \eta^{2l} \dx =
\iom |\partial_{11} \udel|\gam{1}{\frac{\alpha}{2}-\frac{\mu_1}{4}}
\gam{1}{s-\frac{\alpha}{2} + \frac{\mu_1}{4}} \eta^{2l} \dx}\nonumber\\
&\leq & c\Bigg[\iom |\partial_{11} \udel|^2 \gam{1}{\alpha - \frac{\mu_1}{2}} \eta^{2l} \dx
+\iom \gam{1}{2s-\alpha + \frac{\mu_1}{2}} \eta^{2l} \dx\Bigg] \psp .
\end{eqnarray}
Here the first integral on the right-hand side is handled with the help of the inequality \reff{proof 1 3} given in 
Proposition \ref{prop cacc}, hence \reff{proof 1 9} implies
\begin{equation}\label{proof 1 10}
\iom |\partial_{11} \udel| \gam{1}{s} \eta^{2l} \dx \leq
c \Bigg[ \iom |\nabla \eta|^2\eta^{2l-2} \gam{1}{\alpha +1} \dx 
+ \iom \gam{1}{2s - \alpha + \frac{\mu_1}{2}}\eta^{2l} \dx \Bigg] \psp .
\end{equation}
Inserting \reff{proof 1 10} in \reff{proof 1 8} yields
\begin{eqnarray}\label{proof 1 11}
\iom \gam{1}{s+1} \eta^{2l} \dx &=& \iom \gam{1}{s} \eta^{2l} \dx 
+ \iom |\partial_1 \udel|^2 \gam{1}{s} \eta^{2l} \dx\nonumber\\
&\leq& c \Bigg[ \iom \big(\eta^{2l} + |\nabla \eta|^2\eta^{2l-2}\big) \gam{1}{s} \dx 
+ \iom |\nabla \eta|^{2}\eta^{2l-2} \gam{1}{\alpha +1} \dx\nonumber\\
&&+ \iom \gam{1}{2s-\alpha + \frac{\mu_1}{2}} \eta^{2l} \dx \Bigg]\psp .
\end{eqnarray}

We now choose $l$ sufficiently large in order to absorb the first two integrals on the r.h.s.~of
\reff{proof 1 11}: consider numbers $\gamma_1$, $\gamma_2 >0$ and choose $\bar{l}\in \nz$ such that
\[
p := \frac{\gamma_1}{\gamma_2} > 1 \fsp\mbox{and} \fsp p \bar{l} \geq l \psp .
\]
Let $q= p/(p-1)$.
Then we have by Young's inequality for any $\eps > 0$
\begin{eqnarray}\label{proof 1 12}
\iom c(\eta) \eta^{2\bar{l}} \gam{1}{\gamma_2} \dx &=& 
\iom c(\eta) \eta^{2\frac{l}{p}} \eta^{2\bar{l} - 2 \frac{l}{p}} \gam{1}{\gamma_2}\dx\nonumber\\
&\leq & \eps \iom \eta^{2l} \gam{1}{\gamma_1}\dx + c(\eps,\eta)
\iom \eta^{\big(2\bar{l} - 2 \frac{l}{p}\big)q} \dx \psp .
\end{eqnarray}

We apply \reff{proof 1 12} with the choices $\gamma_1 = s+1$ and $\gamma_2 = s$, $\gamma_2 = \alpha +1$,
respectively, recalling $\alpha < s$.\\

Moreover for $\bar{l} := l-1$ we have that $p \bar{l} \geq l$, if $l$ sufficiently large, and \reff{proof 1 9} finally gives
\begin{equation}\label{proof 1 13}
\iom \eta^{2l} \gam{1}{s+1} \dx \leq c \Bigg[1+\iom \eta^{2l} \gam{1}{2s-\alpha +\frac{\mu_1}{2}} \dx\Bigg] \psp .
\end{equation}

We note that
\[
2s - \alpha + \frac{\mu_1}{2} = s + \frac{\hat{\eps}}{2}+\frac{\mu_1}{2} = s + \frac{2+\mu_1}{2}\psp ,
\]
thus with \reff{proof 1 13} we have Proposition  \ref{prop high}.  \qed\\
To finish the proof of Theorem \ref{intro theo 1} we recall $\mu_1 < 2$ and write \reff{proof 1 7} in the form
\begin{equation}\label{proof 1 14}
\iom \eta^{2l} \gam{1}{\gamma_1}\dx \leq c \Bigg[ 1 + \iom \eta^{2l} \gam{1}{\gamma_2}\dx \Bigg]\psp , \fsp
\frac{\gamma_1}{\gamma_2} > 1 \psp .
\end{equation}

Then the same way of absorbing terms as outlined in \reff{proof 1 12} completes the proof of our main theorem. \qed.

\section{Proof of Corollary \ref{intro cor 1}}\label{prcor}

Clearly we may apply the lines of Theorem \ref{intro theo 1} both for $\partial_1 u$ and for $\partial_2 u$
and obtain on account of
\[
\partial_i u_\delta \in L^\chi_{\op{loc}}(\Omega)\fsp\mbox{for $i=1$, $2$, for all $\chi$ and uniform in $\delta$}
\]
for any $\alpha_i \geq 0$, $i=1$, $2$:
\begin{equation}\label{prcor 1}
\iom D^2 f_\delta \big(\nabla \partial_i \udel, \nabla \partial_i \udel\big) \gam{i}{\alpha_i}\eta^{2} \dx \leq c \psp ,
\end{equation}
where $c=c(\eta,\alpha_1,\alpha_2)$ denotes a local constant independent of $\delta$.\\

Given \reff{prcor 1} let us shortly discuss the stress tensor $\sigma$, i.e.~the solution of the
dual variational problem.
In \cite{Bi:2000_1} (see also Section 2.2 of \cite{Bi:1818}) it is shown by elementary arguments from measure theory, 
that the dual problem admits a unique solution and this in turn will give the uniqueness of generalized minimizers up to additive constants
as it will we outlined below.\\

We note that as a general hypothesis of \cite{Bi:2000_1} it is supposed that
\begin{equation}\label{prcor 2}
0 \leq D^2f(\xi)(\eta,\eta) \leq c (1+|\xi|^2)^{-\frac{1}{2}} |\eta|^2 \psp ,
\end{equation}
where the second inequality is not valid in the setting under consideration.\\

However, following the proof of  \cite{Bi:2000_1}, condition \reff{prcor 2} is just needed for showing the uniform
local $W^{1}_2$-regularity of the regularized sequence $\sigma_\delta$ which now immediately follows from
\reff{prcor 1} by choosing $\alpha_1$, $\alpha_2$ sufficiently large. 
As an important consequence we also have Theorem A.9 of \cite{Bi:1818}.\\

Now we claim that for any $\Omega' \Subset \Omega$ and uniformly in $\delta$
\begin{equation}\label{prcor 3}
i)\msp \|\nabla^2 \udel\|_{L^2(\Omega' ;\rz^{2\times 2})} \leq c \psp ,\fsp 
ii)\msp \|\nabla \udel\|_{L^{\infty}(\Omega' ;\rz^2)} \leq c 
\end{equation}
with a constant $c=c(\Omega')$.\\

In fact, we have for arbitrary exponents $\alpha_1$, $\alpha_2 >0$
\begin{eqnarray*}
\lefteqn{\iom \Bigg[ \gam{1}{\alpha_1-\frac{\mu_1}{2}}|\partial_{11} \udel|^2 +
\gam{2}{\alpha_2-\frac{\mu_2}{2}} |\partial_{22} \udel|^2}\\
&+& \Big(\gam{2}{-\frac{\mu_2}{2}}\gam{1}{\alpha_1} + \gam{1}{-\frac{\mu_1}{2}} \gam{2}{\alpha_2}\Big)
|\partial_1\partial_2 \udel|^2\Bigg] \eta^2\dx\\
&&\leq \iom\Bigg[f_1''(\partial_1 \udel) |\partial_{11} \udel|^2 
+ f_2''(\partial_2 \udel) |\partial_{1}\partial_ 2 \udel|^2\Bigg]\gam{1}{\alpha_1}\eta^2 \dx\\
&&\fsp +\iom\Bigg[f_1''(\partial_1 \udel) |\partial_{1}\partial_2 \udel|^2 
+ f_2''(\partial_2 \udel) |\partial_{22} \udel|^2\Bigg]\gam{2}{\alpha_2}\eta^2 \dx\\
&&\leq c \sum_{i=1}^2\iom D^2f_\delta(\nabla \udel)\big(\nabla \partial_i \udel,\nabla \partial_i \udel\big) \gam{i}{\alpha_i}\eta^2 \dx 
\end{eqnarray*}
and by \reff{prcor 1} we obtain the first claim of \reff{prcor 3}. For the second claim we refer, for instance, to
Theorem 5.22 of \cite{Bi:1818}.\\

Given \reff{prcor 3}, we pass to the limit in the differentiated Euler equation for $\udel$ and obtain
\[
\iom D^2f(\nabla u)(\nabla v,\nabla \varphi)\dx = 0 \fsp\mbox{for all}\msp \varphi\in C^1_0(\Omega)
\]
for the function $v=\partial_i u$, $i=1$, $2$.
Observing that the coefficients in this equation are locally bounded and uniformly elliptic, 
a standard reasoning (see, e.g., \cite{GT:1998_1}, Theorem 8.22)
implies H\"older continuity of $v$.\\

Then, by the ``stress-strain relation'' for the particular generalized minimizer and the unique dual solution $\sigma$,
\[
\sigma = \nabla f(\nabla u) \psp ,
\]
we have continuity of $\sigma$ and  $\sigma$ takes values in the set $\op{Im}\nabla f$.\\

Thus, we may apply Theorem A.9 of \cite{Bi:1818} to obtain the uniqueness of generalized minimizers up to an
additive constant.  \hspace*{\fill}\qed\\

\section{Proofs of Corollary \ref{intro cor 2}, Corollary \ref{intro cor 3} and Theorem \ref{intro theo 2}}\label{proof vari}

\emph{Proof of Corollary \ref{intro cor 2}.} We first note that the corollary does not require a splitting structure and that
in the case $-1 < \kappa < 0$ we may take $q:= 2 - \kappa$ as well as
\[
f_\delta(\xi) := \frac{\delta}{q} |\xi|^q + f(\xi)\psp , \fsp \xi =\big(\xi_1,\xi_2\big) \in \rz^2\psp .
\]
Then, with an obvious meaning of  $\udel$, with $\eta\in C^\infty_0(\Omega)$, $0\leq \eta\leq 1$, and 
by letting $\Gamma_\delta := 1+|\nabla \udel|^2$, it holds
for any $\alpha \geq 0$, $l\in \nz$
\begin{equation}\label{cor 2 1}
\iom D^2 f_\delta (\nabla \udel)  \Big(\partial_i \nabla \udel ,\nabla\big[\eta^{2l} \partial_i \udel \Gamma_\delta^{\alpha}\big]\Big) \dx = 0 \psp ,
\end{equation}
where now the sum is taken w.r.t. $i=1$, $2$.\\

From \reff{cor 2 1} we derive as counterpart to \reff{proof 1 3} of Proposition \ref{prop cacc}
\begin{equation}\label{cor 2 2}
\iom \eta^{2l} |\nabla^2 \udel|^2 \Gamma_\delta^{\alpha - \frac{\mu}{2}}\dx \leq
c \iom |\nabla \eta|^{2}\eta^{2l-2} \Gamma_\delta^{\alpha + 1 - \frac{\kappa}{2}}\dx \psp .
\end{equation}
We note that $q$ was defined in such a way that we do not have to consider an extra $\delta$-part on the the right-hand side
of \reff{cor 2 2}.\\

Next \reff{proof 1 11} has to be replaced by
\begin{eqnarray}\label{cor 2 3}
\iom \Gamma_\delta^{s+1} \eta^{2l} \dx
&\leq& c \Bigg[ \iom \big(\eta^{2l} + |\nabla \eta|^2\eta^{2l-2}\big) \Gamma_\delta^{s} \dx \nonumber\\
&&+ \iom |\nabla \eta|^{2}\eta^{2l-2} \Gamma_\delta^{\alpha +1-\frac{\kappa}{2}} \dx\nonumber\\
&&+ \iom \Gamma_\delta^{2s-\alpha + \frac{\mu_1}{2}} \eta^{2l} \dx \Bigg]\psp .
\end{eqnarray}

Then we may proceed as in the proof of Proposition \ref{prop high} provided that the exponents $\alpha$ and $s$ satisfy the conditions
\begin{eqnarray}\label{cor 2 4}
\alpha + 1 - \frac{\kappa}{2} < s+1 &\Leftrightarrow & \alpha < s +\frac{\kappa}{2}\\[2ex]
\label{cor 2 5}
2s - \alpha + \frac{\mu}{2} < s+1 & \Leftrightarrow & s < \alpha + \frac{2-\mu_1}{2}\psp .
\end{eqnarray}
But \reff{cor 2 4} and \reff{cor 2 5} are consequences of  \reff{intro 17} which proves the corollary. \hspace*{\fill}\qed\\

\emph{Proof of Corollary \ref{intro cor 3}.} For instance we now may consider $q:= 2 + \varkappa$ as well as
\[
f_\delta(\xi) := \frac{\delta}{q} |\xi_1|^q + \frac{\delta}{2}|\xi|^2 + f(\xi)\psp , \fsp \xi =\big(\xi_1,\xi_2\big) \in \rz^2\psp .
\]

Going through the proof of Theorem \ref{intro theo 1} and adapting the exponents according to the proof of
\reff{intro cor 2} we obtain \reff{intro cor 3} by recalling the boundedness of $f_2''$. \hspace*{\fill}\qed\\ 

\emph{Proof of Theorem \ref{intro theo 2}.}
With obvious changes in notation we now use as regularizing energy density:
\[
f_{1,\delta} (t) :=\frac{\delta}{2} t^2 + f_1(t) \psp ,\fsp
f_{2,\delta}(t) := \frac{1}{\gamma +2} |t|^{\gamma +2} \psp .
\]

Then we are going to establish a variant of Proposition \ref{prop cacc} adapted to the hypothesis  \reff{vari 4}. 
As the main new feature, we consider negative exponents in the Caccioppoli-inequality. More precisely, we have:

\begin{proposition}\label{prop vari 1}
Fix $l\in \nz$ and suppose that $\eta \in C^\infty_0(\Omega)$, $0 \leq \eta \leq 1$. Then, given the assumptions of Theorem \ref{intro theo 2},
the inequality
\begin{eqnarray}\label{proof vari 1}
\lefteqn{\iom D^2 f_\delta(\nabla \udel)\big(\nabla \partial_1 \udel, \nabla \partial_1 \udel\big) \eta^{2l} 
\gam{1}{\alpha} \dx}\nonumber \\ 
&& \leq c \iom D^2f_\delta(\nabla \udel) (\nabla \eta,\nabla \eta)\eta^{2l-2} \gam{1}{\alpha +1}\ \dx
\end{eqnarray}
holds for any $\alpha > - 1/2$, which in particular implies (again for all $\alpha > -1/2$)
\begin{eqnarray}\label{proof vari 2}
\lefteqn{\iom \eta^{2l} \gam{1}{\alpha - \frac{\mu_1}{2}} |\partial_{11}\udel|^2 \dx}\nonumber\\
&\leq& c \Bigg[ 1+ \iom |\nabla \eta|^2 \eta^{2l-2}\gam{1}{\alpha +1} \dx  
+ \iom |\nabla \eta|^2 \eta^{2l-2}\gam{1}{\frac{\alpha+1}{1-\gamma}}\dx\Bigg] \psp .
\end{eqnarray}
As usual $c=c(l,\alpha)$ denotes a uniform constant independent of $\delta$.
\end{proposition}

\emph{Proof of Proposition \ref{prop vari 1}.} Suppose w.l.o.g.~that $-1/2 < \alpha \leq 0$. 
Exactly as outlined in the proof of Proposition \ref{prop cacc} we have \reff{proof 1 4}
together with the inequality 
\begin{eqnarray*}
|S_1| &=&  2 |\alpha| \iom D^2f_\delta (\nabla \udel) \big(\nabla \partial_1 \udel, \nabla \partial_1 \udel\big)
|\partial_1 \udel|^2\gam{1}{\alpha-1}\eta^{2l}\dx\\
&\leq & 2  |\alpha| \iom D^2f_\delta (\nabla \udel) \big(\nabla \partial_1 \udel, \nabla \partial_1 \udel\big)\gam{1}{\alpha}\eta^{2l}\dx \psp .
\end{eqnarray*}
On account of $2 |\alpha| < 1$ we may absorb $|S_1|$ on the left-hand side of \reff{proof 1 4} with the result
\[
\iom D^2f_\delta (\nabla \udel) \big(\nabla \partial_1 \udel, \nabla \partial_1 \udel\big) \eta^{2l}\gam{1}{\alpha}\dx \leq c |S_2| \psp .
\]
Here $|S_2|$ may be estimated in the same manner as in the proof of Propostion \ref{prop cacc} leading to \reff{proof vari 1}.\\

A variant of inequality \reff{proof 1 6} is based on the assumptions \reff{vari 3} and \reff{vari 4}. We have
\begin{eqnarray}
\lefteqn{\iom D^2 f_\delta(\nabla \udel)\big(\nabla \eta,\nabla \eta\big) \eta^{2l-2} \gam{1}{\alpha+1}\dx}\nonumber\\
&\leq & c \Bigg[ \iom |\nabla \eta|^2 \eta^{2l-2} \gam{1}{\alpha+1}\dx 
+ \iom |\nabla \eta|^2 \eta^{2l-2}\gam{2}{\frac{\gamma}{2}} \gam{1}{\alpha+1}\dx\Bigg]\psp .\label{proof vari 3}
\end{eqnarray}
Letting
\[
1 < p_1= \frac{1}{\gamma}\psp , \fsp p_2 = \frac{1}{1-\gamma} \psp ,
\]
we apply Young's inequality with the result
\begin{eqnarray}
\iom |\nabla \eta|^2 \eta^{2l-2}\gam{2}{\frac{\gamma}{2}} \gam{1}{\alpha+1}\dx
&=& \iom \Big[|\nabla \eta|^2 \eta^{2l-2}\Big]^{\frac{1}{p_1}} \gam{2}{\frac{\gamma}{2}} \Big[|\nabla \eta|^2 \eta^{2l-2}\Big]^{\frac{1}{p_2}} \gam{1}{\alpha+1}\dx\nonumber\\
&\leq & c \Bigg[ 1+ \iom |\nabla \eta|^2 \eta^{2l-2} \gam{1}{\frac{\alpha+1}{1-\gamma}} \dx \Bigg] \psp .\label{proof vari 4}
\end{eqnarray}
 With \reff{proof vari 3} and \reff{proof vari 4} the proof of Proposition \ref{prop vari 1} is completed. \qed.\\

\begin{proposition}\label{prop vari 2}
Let the hypotheses of Theorem \ref{intro theo 2} hold, in particular we have \reff{vari 5}. Consider
real numbers $\tau_s$, $\tau_\alpha >0$ such that
\begin{equation}\label{proof vari 5}
s := -\frac{1}{2}+\tau_s \psp , \fsp \alpha := -\frac{1}{2} + \tau_\alpha  \psp . 
\end{equation}
Then, for $l \in \nz$ and with a local constant $c(\eta,l)$ we have
\begin{eqnarray}\label{proof vari 6}
\iom \eta^{2l} \gam{1}{s+1} \dx &\leq& c \Bigg[1+
\iom |\nabla \eta|^2 \eta^{2l-2}\gam{1}{\alpha +1} \dx\nonumber\\  
&&+ \iom |\nabla \eta|^2 \eta^{2l-2}\gam{1}{\frac{\alpha+1}{1-\gamma}} \dx\nonumber\\ 
&&+\iom \eta^{2l} \gam{1}{2s-\alpha+\frac{\mu_1}{2}} \dx\Bigg] \psp .
\end{eqnarray}
\end{proposition}

\emph{Proof of Proposition \ref{prop vari 2}.} Exactly as in the proof of Proposition \ref{prop high} we obtain
inequalities \reff{proof 1 8} and \reff{proof 1 9}.\\

Now \reff{proof 1 10} has to be replaced by
\begin{eqnarray}\label{proof vari 7}
\iom |\partial_{11} \udel| \gam{1}{s} \eta^{2l} \dx &\leq&
c \Bigg[ 1+ \iom |\nabla \eta|^2 \eta^{2l-2}\gam{1}{\alpha +1} \dx  \nonumber\\
&&+ \iom |\nabla \eta|^2 \eta^{2l-2}\gam{1}{\frac{\alpha+1}{1-\gamma}}\dx\nonumber\\
&&+ \iom \gam{1}{2s - \alpha + \frac{\mu_1}{2}}\eta^{2l} \dx \Bigg] \psp .
\end{eqnarray}
Inserting \reff{proof vari 7} in \reff{proof 1 8} yields the claim of the proposition.
\hspace*{\fill}\qed\\

In order to finish the proof of Theorem \ref{intro theo 2} with similar arguments as applied in
Theorem  \ref{intro theo 1}, we have to verify the conditions
\begin{eqnarray}\label{proof vari 8}
\alpha + 1 < s+1 & \Leftrightarrow & \tau_\alpha < \tau_s \psp ,\\[2ex]
\frac{\alpha +1}{1-\gamma} < s + 1 &\Leftrightarrow&   \gamma < \frac{2 (\tau_s-\tau_\alpha)}{1+2 \tau_s} \psp ,\label{proof vari 9}\\[2ex]
2s- \alpha + \frac{\mu_1}{2} < s+ 1 & \Leftrightarrow & \tau_s - \tau_\alpha < 1 - \frac{\mu_1}{2} \psp .\label{proof vari 10}
\end{eqnarray}
But this can be done as follows: for any
\begin{equation}\label{proof vari 11}
\tau_s <  1 - \frac{\mu_1}{2}
\end{equation}
we choose $\tau_\alpha > 0$ sufficiently small such that \reff{proof vari 10} holds. Since we have \reff{vari 5},
we may also increase $\tau_s$, still satisfying \reff{proof vari 11}, such that \reff{proof vari 9} and
trivially \reff{proof vari 8} are satisfied. This completes the proof of Theorem \ref{intro theo 2}.\hspace*{\fill}\qed\\

\section{Appendix}\label{app}

\subsection{Proof of  Lemma \ref{intro lem 1}}\label{app 1}
Arguing by contradiction we assume the validity of the second inequality in \reff{intro 9} with some exponent
$\kappa > 1$.\\

W.l.o.g.~we assume $n=2$ since otherwise we can replace $f$ by $\tilde{f}$: $\rz^2 \to \rz$,
$\tilde{f}(p_1,p_2) := f(p_1,p_2, 0 ,\dots ,0)$ observing that $\tilde{f}$ satisfies \reff{intro 3},
\reff{intro 9} for $\xi$, $\eta\in \rz^2$ with the same exponents $\mu$ and $\kappa$.\\

Consider an increasing sequence of numbers $c_k >0$ such that
\begin{equation}\label{app 1 1}
\lim_{k\to \infty} c_k = \infty
\end{equation}
and let
\[
A_k := \big[ f \leq c_k\big] := \big\{ p \in \rz^2: f(p) \leq c_k\big\} \psp .
\]

Condition \reff{intro 3} implies the boundedness of each set $A_k$, moreover, we have strict
convexity of $A_k$ on account of \reff{intro 9}.\\

Let $\gamma_k$ denote a parametrization by arc length of the closed convex curve
$\partial A_k = \big[f = c_k]$ -- in fact, we just need a parametrization inside 
a small neighborhood of the point $p_k$ considered in \reff{app 1 6} below.
For each $k$ it holds $c_k= f\big(\gamma_k(t)\big)$,
hence
\begin{eqnarray}
0&=& \gamma_k'(t) \cdot \gamma_k ''(t) \psp , \label{app 1 2}\\[1ex]
0&=& \gamma_k'(t) \cdot Df\big(\gamma_k(t)\big) \psp , \label{app 1 3}
\end{eqnarray}
and \reff{app 1 3} implies by taking the derivative
\begin{equation}\label{app 1 4}
0 = D^2f\big(\gamma_k(t)\big) \big(\gamma_k'(t),\gamma_k'(t)\big) + Df \big(\gamma_k(t)\big) \cdot
\gamma_k''(t) \psp .
\end{equation}

From \reff{app 1 2} and \reff{app 1 3} we deduce that $Df\big(\gamma_k(t)\big)$ and $\gamma_k''(t)$ are
proportional for each $t$, therefore \reff{app 1 4} implies
\[
\big|\gamma_k''(t)\big| \psp \big| Df\big(\gamma_k(t)\big) \big| =
\big| \gamma_k''(t) \cdot Df_k\big(\gamma_k(t)\big)\big| =
D^2 f\big(\gamma_k(t)\big) \big(\gamma_k'(t),\gamma_k'(t)\big) \psp .
\]

Combing this equation with \reff{intro 9} we find
\begin{equation}\label{app 1 5}
\big|\gamma_k''(t)\big|\psp \big|Df\big(\gamma_k(t)\big)\big| \leq 
c_2 \big( 1+|\gamma_k(t)|\big)^{-\kappa} \psp .
\end{equation}

Recalling the boundedness of $A_k$ and assuming as usual $f(0)=0$ we find a radius
$r_k$ such that $\overline{A}_k \subset \overline{B_{r_k}}(0)$ and in addition with the property
that $\partial A_k$ and $\partial B_{r_k}(0)$ have at least one common point
$p_k = \gamma_k(t_k)$. Then it holds
\begin{equation}\label{app 1 6}
\big|\gamma''_k(t_k)\big|\geq \frac{1}{r_k} \psp ,
\end{equation}
moreover, we have $\big|\gamma_k(t_k)\big| = r_k$. Combining \reff{app 1 5} and \reff{app 1 6}
we obtain
\begin{equation}\label{app 1 7}
\frac{1}{\big| \gamma_k(t_k)\big|} \big|Df\big(\gamma_k(t_k)\big)\big| 
\leq c_2 \big(1+|\gamma_k(t_k)|\big)^{-\kappa} \psp .
\end{equation}
The convexity of $f$ implies
\begin{equation}\label{app 1 8}
Df\big(\gamma_k(t_k)\big) \cdot \gamma_k(t_k) \geq b_1 \big|\gamma_k(t_k)\big| - b_2
\end{equation}
with suitable real numbers $b_1 >0$, $b_2 \geq 0$ independent of $k$.
Recalling $f\big(\gamma_k(t)\big) = c_k$ as well as our assumption \reff{app 1 1} concerning the
sequence $(c_k)$, condition \reff{intro 3} implies $\big|\gamma_k(t_k)\big| \to \infty$, hence
\reff{app 1 8} shows for $k \gg 1$ the validity of
\begin{equation}\label{app 1 9}
Df\big(\gamma_k(t_k)\big)\cdot \gamma_k(t_k) \geq \frac{1}{2} b_1 \big|\gamma_k(t_k)\big| \psp .
\end{equation}
From \reff{app 1 7} and \reff{app 1 9} it finally follows
\[
\frac{1}{\big|\gamma_k(t_k)\big|} \leq \frac{2}{b_1} 
\frac{\big|Df\big(\gamma_k(t_k)\big)\cdot \gamma_k(t_k)\big|}{\big|\gamma_k(t_k)\big|^2} 
\leq \frac{2}{b_1} c_2 \big(1+|\gamma_k (t_k)|\big)^{-\kappa} \psp ,
\]
or equivalently
\[
\frac{b_1}{2 c_2} \leq \big| \gamma_k(t_k)\big| \big(1+|\gamma_k(t_k)|\big)^{-\kappa} \psp .
\]
Note that the r.h.s.~vanishes as $k\to \infty$ in case $\kappa > 1$ leading to a contradiction and
proving the lemma.\hfill \qed


\subsection{Examples}\label{app 2}
Even in the one-dimensional setting there is a fundamental difference to the discussion of examples with superlinear growth, where
a natural scale is given by $h_p(t) = (1+t)^p$, $t\geq 0$, $p>1$, with growth rate $p$ and satisfying a balancing condition \reff{intro 8}
with exponent $p-2 = - \mu$ in both directions.\\

In the case of linear growth problems we have, for instance, the functions $\Phi_\mu$ satisfying \reff{intro 8}. However,
the growth of $\Phi_\mu$ is linear, hence not depending on $\mu$.\\

In this second part of the appendix we are going to sketch some examples which show that the prototypes
$\Phi_\mu$ leading to energy densties with linear growth are far apart from being universal representatives.\\

The construction of examples $h$: $\rz^+_0 \to \rz^+_0$ of class $C^2$ such that $\hat{h}(t)=h(|t|)$ fits into the above discussion
is limited by the two constraints
\begin{itemize}
\item $h''(t) > 0$ for all $t\in \rz^+_0$,
\item $\dis \int_0^\infty h''(\tau) \D \tau < \infty$.
\end{itemize}  

We just mention three examples which are interesting from different points of view.
\begin{enumerate}
\item Let us start by considering the family $\Phi_\mu/(\mu -1)$ in limit case $\mu =1$. In this limit case
we have
\[
\lim_{\mu  \downarrow 1} \frac{\Phi_\mu(t)}{\mu-1} = t \ln(1+t) + \ln(1+t) - t \psp ,
\]
i.e.~the limit case corresponds to an energy density of nearly linear growth. We now present an easy example 
given in an explicit analytic form such that
\[
h''(t) \leq c (1+t)^{-1}\fsp\mbox{for all}\msp t\in \rz^+_0
\]
with optimal exponent $-1$ on the right-hand side and such that $f$ is of linear growth. We let 
\begin{eqnarray*}
h(t) &=& \int_0^t h'(\tau)\D \tau\psp ,\fsp h'(t) = 1- (1+t)^{1-\mu(t)} \psp ,\\
\mu(t) &=& 1 + \frac{1}{\ln(1+t)}\psp , \fsp t \gg 1\psp .
\end{eqnarray*}

Then we have the desired estimate
\[
h''(t) = 2 (1+t)^{-1 -\frac{1}{\ln(1+t)}} \frac{1}{\ln(1+t)}\psp , \fsp t \gg 1 \psp .
\]

\item The second example is of linear growth with a lower bound for the second derivative which is even worse
than involving an exponent $-\mu$. We consider the elementary function
\[
h(t) = \ln\big( 1+ e^t\big) \fsp\mbox{for all}\msp t\in \rz^+_0 \psp .
\]

An easy calculation shows for allt $t\geq 0$
\[
h'(t) = \frac{e^t}{1+e^{t}}\psp ,\fsp h''(t) = \frac{1}{4} \frac{1}{\cosh^2(t/2)}\psp .
\]

\item In view of the above given line for the construction of examples we may also have a countable union of atoms.
This example corresponds to an approximation of a convex piecewise affine continuous function of linear growth.\\

For $i \in \nz$ and $\sigma_i  > 0$ we consider
\[
h''(t) = \sum_{i=1}^{\infty} e^{-\frac{|t-i|^2}{\sigma_i^2}}\psp ,\fsp t\geq 0 \psp ,
\]
and for $\sigma_i$ sufficiently small one obtains
\[
\int_{0}^\infty h''(t) \dt \leq c \sum_{i=1}^\infty \sigma_i  \leq c\psp .
\]
\end{enumerate}



\begin{tabular}{ll}
Michael Bildhauer&bibi@math.uni-sb.de\\
Martin Fuchs&fuchs@math.uni-sb.de\\[5ex]
Department of Mathematics&\\
Saarland University&\\
P.O.~Box 15 11 50&\\
66041 Saarbr\"ucken&\\ 
Germany&
\end{tabular}

\end{document}